\documentclass[11pt]{article}

\usepackage[margin=1in]{geometry}
\usepackage{amsmath,amsfonts,amsthm,latexsym,xspace}
\usepackage{graphicx}
\usepackage{psfrag}
\newtheorem{thm}{Theorem}[section]
\newtheorem{cor}[thm]{Corollary}
\newtheorem{lem}[thm]{Lemma}

\theoremstyle{definition}

\theoremstyle{remark}
\newtheorem{rem}[thm]{Remark}
\newcommand{\set}[1]{\left\{#1\right\}}

\newcommand{\eps}{\varepsilon}

\newcommand{\e}{\mathbf{E}}
\newcommand{\M}{\mathcal{M}}
\newcommand{\C}{\mathcal{C}}
\newcommand{\CE}{\mathcal{E}}
\newcommand{\CG}{\mathcal{G}}
\newcommand{\CH}{\mathcal{H}}
\newcommand{\CI}{\mathcal{I}}
\newcommand{\CN}{\mathcal{N}}
\newcommand{\CV}{\mathcal{V}}
\newcommand{\binpower}[2]{#1\rule{0in}{2.2ex}^{{\binom{#2}{2}}}}


\newcommand{\hide}[1]{}  % replace {} by {#1} to reveal

\title{Path Coupling Using Stopping Times and Counting\\Independent Sets
and Colourings in Hypergraphs}
\author{Magnus Bordewich\thanks{School of
Computing, University of Leeds, Leeds LS2 9JT, UK. Email:
\texttt{\{dyer,magnusb\}@comp.leeds.ac.uk}.},\ \, Martin Dyer${}^*$ and Marek
Karpinski\thanks{Dept. of Computer Science, University of Bonn, 53117 Bonn,
Germany. Email: \texttt{marek@cs.uni-bonn.de}.}}
\date{April 2, 2005}
\begin{document}
\maketitle

\begin{abstract}
We analyse the mixing time of Markov chains using
path coupling with stopping times. We apply this approach to two hypergraph
problems. We show that the Glauber dynamics for independent sets in a
hypergraph mixes rapidly as long as the maximum degree $\Delta$ of a vertex
and the minimum size $m$ of an edge satisfy $m\geq 2\Delta+1$. We also show
that the Glauber dynamics for proper $q$-colourings of a hypergraph mixes
rapidly if $m\geq 4$ and $q > \Delta$, and if $m=3$ and $q\geq1.65\Delta$.
We give related results on the hardness of exact and approximate counting
for both problems.
\end{abstract}

\section{Introduction}

We develop a  new approach to using stopping times in conjunction with path
coupling to bound the convergence of time of Markov chains. Our main
interest is in applying these results to randomised approximate counting.
For an introduction, see~\cite{J03}. To illustrate our methods, we consider
approximation of the numbers of independent sets and $q$-colourings in
hypergraphs with upper-bounded degree, and lower-bounded edge size. These
problems in hypergraphs are of interest in their own right but, while
approximate optimisation has received
attention~\cite{DRS02,DGKR03,HL98,KNS01}, there has been surprisingly
little work on approximate counting.

Our results are achieved by considering, in the path coupling setting, the
stopping time at which the distance between two coupled chains first
changes.  The first application of stopping times to path coupling was by
Dyer, Goldberg, Greenhill, Jerrum and Mitzenmacher~\cite{DGGJM01}. Their
analysis was later improved by Hayes and Vigoda~\cite{HV04}, using a method
closely related to that developed in this paper. Theorem~\ref{stopping},
the main technical result of the paper, shows that if the expected distance
between the two chains has decreased at this stopping time, then the chain
mixes rapidly. This also follows from~\cite[Corollary~4]{HV04}.  However we
give a simpler proof than that of~\cite{HV04}, and our
Theorem~\ref{stopping} will usually give a moderate improvement in the
bound on mixing time in comparison with~\cite[Corollary~4]{HV04}. See
Remark~\ref{rem20} below.

The problem of approximately counting independent sets in graphs has been
widely studied, see for example~\cite{DFJ99,DG00a,LV99,M01,V01}, but the
only previous work on the approximate counting of independent sets in
\emph{hypergraphs} seems to that of Dyer and Greenhill~\cite{DG00a}. They
showed rapid mixing to the uniform distribution of a simple Markov chain on
independent sets in a hypergraph with maximum degree 3 and maximum edge
size 3. However, this was the only interesting case resolved. Their results
imply rapid mixing only for $m\leq \Delta/(\Delta-2)$, which gives $m\leq
3$ when $\Delta=3$ and $m\leq 2$ when $\Delta\geq 4$. In
Theorem~\ref{indsets} we prove rapid mixing of the \emph{Glauber dynamics}
for any hypergraph such that $m\geq 2\Delta+1$, where $m$ is the smallest
edge size and $\Delta$ is the maximum degree. This is a marked improvement
for large $m$. More generally, we consider the \emph{hardcore distribution}
on independent sets with \emph{fugacity} $\lambda$. (See, for
example,~\cite{DG00a,LV99,V01}.) In~\cite{DG00a}, it is proved that rapid
mixing occurs if $\lambda\leq m/((m-1)\Delta-m)$. Here we improve this
considerably for larger values of $m$, to $\lambda\leq (m-1)/2\Delta$. We
also give proofs that computing the number of independent sets in
hypergraphs is \#P-complete except in trivial cases, and that there can be
no approximation for the number of independent sets in a hypergraphs if the
minimum edge size is at most logarithmic in $\Delta$. It may be noted that
our upper and lower bounds are exponentially different. We have no strong
belief that either is close to the threshold at which approximate counting
is possible, if such a threshold exists.

Counting $q$-colourings of hypergraphs was considered by Bubley~\cite{B01},
who showed that the Glauber dynamics was rapidly mixing if $q\geq 2\Delta$,
generalising a result of Jerrum~\cite{J95} and Salas and Sokal~\cite{SS97}
for graphs.  Much work has been done on improving this result for graph
colourings, see~\cite{DFHV04} and its references, but little attention
appears to have been given to the hypergraph case. Here we prove rapid
mixing of Glauber dynamics for proper colourings of hypergraphs if  $m\geq
4$, $q>\Delta$, and if $m=3$, $q\geq1.65\Delta$. For a precise statement of
our result see Theorem~\ref{colouring}. Again we give proofs that computing
the number of colourings in hypergraphs is \#P-complete except in trivial
cases, and that there can be no approximation for the number of colourings
of hypergraphs if $q\leq (1-1/m)\Delta^{1/(m-1)}$. Again, there is a
considerable discrepancy between the upper and lower bounds for large $m$.

The paper is organised as follows. Section~\ref{sec:intuit} gives an
intuitive motivation for the stopping time approach of the paper.
Section~\ref{sec:stopping} contains the full description and proof of
Theorem~\ref{stopping} for path coupling with stopping times. We apply this
to hypergraph independent sets in Section~\ref{sec:indsets}.
Section~\ref{sec:hard} contains the hardness proofs.
Section~\ref{sec:colour} contains analysis of the Glauber dynamics for
hypergraph colouring. Finally, Section~\ref{sec:hardcol} contains the
hardness results for counting colourings in hypergraphs.

\subsection{Intuition}\label{sec:intuit}
Let $\CH=(\CV,\CE)$ be a hypergraph of maximum degree $\Delta$ and minimum
edge size $m$. A subset $S\subseteq \CV$ of the vertices is
\emph{independent} if no edge is a subset of $S$. Let $\Omega(\CH)$ be the
set of all independent sets of $\CH$. Let $\lambda$ be the \emph{fugacity},
which weights independent sets. (See~\cite{DG00a}.) The most important case
is $\lambda=1$, which weights all independent sets equally and gives rise
to the uniform distribution on all independent sets. We define the Markov
chain $\M(\CH)$ with state space $\Omega(\CH)$ by the following transition
process (\emph{Glauber dynamics}). If the state of $\M$ at time $t$ is
$X_t$, the state at $t+1$ is determined by the following procedure.
\begin{enumerate}
\item Select a vertex $v\in \CV$ uniformly at random, \item
\begin{enumerate}
    \item if $v\in X_t$ let $X_{t+1}=X_t\backslash \{v\}$ with probability
    $1/(1+\lambda)$,
    \item if $v\not\in X_t$ and $X_t \cup \{v\} $ is independent, let
    $X_{t+1}=X_t\cup \{v\}$ with probability $\lambda/(1+\lambda)$,
    \item otherwise let $X_{t+1}=X_t$.
\end{enumerate}
\end{enumerate}
This chain is easily shown to be ergodic with stationary probability
proportional to $\lambda^{|I|}$ for each independent set $I\subseteq \CV$.
In particular, $\lambda=1$ gives the uniform distribution. The natural
coupling for this chain is the ``identity'' coupling, the same transition
is attempted in both copies of the chain. If we try to apply standard path
coupling to this chain, we immediately run into difficulties. Consider two
chains $X_t$ and $Y_t$ such that $Y_t=X_t\cup\set{w}$, where $w\notin X_t$
(the \emph{change vertex}) is of degree $\Delta$. An edge $e\in\CE$ is
\emph{critical} in $Y_t$ if it has only one vertex $z\in\CV$ which is not
in $Y_t$, and we call $z$ \emph{critical for $e$}. If each of the edges
through $w$ is critical for $Y_t$, then there are $\Delta$ choices of $v$
in the transition which can be added in $X_t$ but not in $Y_t$. Thus, if
$\lambda=1$, the change in the expected Hamming distance between $X_t$ and
$Y_t$ after one step could be as high as $\frac{\Delta}{2n}-\frac1n$. Thus
we obtain rapid mixing only in the case $\Delta=2$. This case has some
intrinsic interest, since the complement of an independent set corresponds,
under hypergraph duality, to an \emph{edge cover}~\cite{GJ79} in a graph.
Thus we may uniformly generate edge covers, but the scope for unmodified
path coupling is obviously severely limited.

The insight on which this paper is based is as follows. Although in one
step it could be more likely that a \emph{bad vertex} (increasing Hamming
distance) is chosen than a \emph{good vertex} (decreasing Hamming
distance), it is even more likely that one of the other vertices in an edge
containing $w$ is chosen and removed from the independent set. Once the
edge has two unoccupied vertices other than $w$, then any vertex in that
edge can be added in both chains. This observation enables us to show that,
if $T$ is defined to be the stopping time at which the distance between
$X_t$ and $Y_t$ first changes, the expected distance between $X_T$ and
$Y_T$ will be less than 1. Theorem~\ref{stopping} below shows that under
these circumstances path coupling can easily be adapted to prove rapid
mixing.

Having established this general result, we use it to prove that $\M(\CH)$
is rapidly mixing for hypergraphs with $m\geq 2\lambda\Delta+1$. Note that,
though all the results in this paper will be proved for uniform hypergraphs
of edge size $m$, they carry through trivially for hypergraphs of minimum
edge size $m$.

\section{Path coupling using a stopping time}\label{sec:stopping}

First we prove the main result discussed above.

\begin{thm}\label{stopping}
Let $\M$ be a Markov chain on state space $\Omega$. Let $\mathrm{d}$ be an integer
valued metric on $\Omega\times \Omega$, and let $(X_t,Y_t)$ be a
path coupling for $\M$, where $S$ is the set of pairs of states $(X,Y)$ such that $\mathrm{d}(X,Y)=1$.
For any initial states $(X_0,Y_0)\in S$ 
let $T$ be the stopping time given by the minimum $t$ such that
$\mathrm{d}(X_t,Y_t)\neq1$. Suppose, for some $p>0$, that
\begin{enumerate}
\item $\Pr(T=t\, |\, T\geq t)\geq p$, independently for each $t$,%
\item $\e[\mathrm{d}(X_T,Y_T)]\leq \alpha < 1$.
\end{enumerate}
Then $\M$ mixes rapidly. In particular the mixing time
$\tau(\eps)$ of $\M$ satisfies
\[\tau(\eps)\ \leq\
\frac{1}{p}\,\frac{3}{1-\alpha} \ln (e
D_2)\ln\Big(\frac{2D_1}{\eps(1-\alpha)}\Big),\]%
where $D_1=\max\{\mathrm{d}(X,Y):X,Y\in\Omega\}$ and
$D_2=\max\{\mathrm{d}(X_T,Y_T):X_0,Y_0\in\Omega,\,
\mathrm{d}(X_0,Y_0)=1\}$.
\end{thm}

\begin{proof}
Consider the following game. In each round a gambler either wins \pounds$
1$, loses some amount \pounds$(l-1)$ or continues to the next round. If he
loses \pounds$(l-1)$ in a game, he starts $l$ separate (but possibly
dependent) games simultaneously in an effort to win back his money. If he
has several games going and loses one at a certain time, he starts $l$ more
games, while continuing with the others that did not conclude. We know that
the probability he finishes a game in a given step is at least $p$, and the
expected winnings in each game is at most $1-\alpha$. The question is: does
his return have positive expectation at any fixed time\,? We will show that
it does. But first a justification for our interest in this game.

Each game represents a single step on the path between two states
of the coupled Markov chain. We start with $X_0$ and $Y_0$
differing at a single vertex. The first game is won if the first
time the distance between the coupled chains changes is in
convergence. The game is lost if the distance increases to $l$. At
that point we consider the distance $l$ path $X_t$ to $Y_t$, and
the $l$ games played represent the $l$ steps in the path. Although
these games are clearly dependent, they each satisfy the
conditions given. The gambler's return at time $t$ is one minus the
length of the path at time $t$, so a positive expected return
corresponds to an expected path length less than one. We will show that the expected path length is sufficiently small to ensure coupling.

First note that the gambler's return at time $t$ is one minus the
number of games active at time $t$. For the initial game we define
the \emph{level} to be zero, for any other possible game we define
the level to be one greater than the level of the game whose loss
precipitated it. We define the random variables $M_k, l_{jk}$ and
$I_{jk(t)}$ as follows. $M_k$ is the number of games at level $k$
that are played, $l_{jk},$ for $j=1\ldots M_k$, is the number of
games in level $k+1$ which are started as a result of the outcome
of game $j$ in level $k$, and $I_{jk}(t)$ is an indicator function
which takes the value 1 if game $j$ in level $k$ is active at time
$t$, and 0 otherwise. Let $N(t)$ be the number of games active at
time $t$. Then, by linearity of expectations,
\begin{equation}
\label{expsum} \e[N(t)]=\sum_{k=0}^\infty \e \left[ \sum_{j=1}^{M_k}
I_{jk}(t)\right].
\end{equation}

We will bound this sum in two parts, splitting it at a point $k=K$
to be determined. For $k\leq K$ we observe that $M_k\leq D_2^k$.
Since $\Pr (I_{jk}(t)=1)$ is at most the probability that exactly
$k-1$ games of a sequence are complete at time $t$, regardless of
outcome, we have
\begin{align*}
\e \left[ \sum_{j=1}^{M_k} I_{jk}(t)\right]\ \leq\ D_2^k \max_j
\e[I_{jk}(t)]\ \leq\ D_2^k \Pr(\textrm{exactly $k-1$ games complete by time
$t$}).
\end{align*}
So that
\begin{align}
\sum_{k=0}^K \e \left[ \sum_{j=1}^{M_k} I_{jk}(t)\right]&\leq\
\sum_{k=0}^K D_2^k \Pr(\textrm{exactly $k-1$ games complete by time $t$})\notag \\
&\leq\ {D_2}^K\Pr(\textrm{at most $K$ games complete by $t$}).\label{front}
\end{align}
On the other hand, for $k>K$ we observe that
\begin{align*}
\e \bigg[ \sum_{j=1}^{M_k} I_{jk}(t)\bigg]&\leq\ \e[M_k]\ =\
\e_{M_{k-1}}\big[\e[M_k|M_{k-1}]\big]\ =\
\e_{M_{k-1}}\big[\e[\sum_{j=1}^{M_{k-1}}l_{jk-1}|M_{k-1}]\big]
\end{align*}
Since $\e[l_{jk-1}]\leq\alpha$ for any starting conditions, we may
apply this bound even when conditioning on $M_{k-1}$. So
\begin{align}
\e \left[ \sum_{j=1}^{M_k} I_{jk}(t)\right]
%&\leq\ \e[\e[\sum_{j=1}^{M_{k-1}}l_{jk-1}|M_{k-1}]]\notag \\
%&\leq\ \e[\e[\sum_{j=1}^{M_{k-1}}\alpha|M_{k-1}]]\notag \\
\ \leq\ \e[\alpha M_{k-1}]\ \leq\ \alpha^k,\label{tail}
\end{align}
using linearity of expectation, induction and $\e[M_1]\leq\alpha$. Putting
\eqref{front} and \eqref{tail} together we get
\begin{align} \e[N(t)]&\leq\ {D_2}^K\Pr(\textrm{at most $K$ games complete by $t$})
\ +\ \sum_{k=K+1}^\infty \alpha^k\notag \\
&=\ {D_2}^K\Pr(\textrm{at most $K$ games complete by $t$}) \ +\
\frac{\alpha^{K+1}}{1-\alpha}.
\end{align}

We now set $K =\lfloor(\ln \alpha)^{-1}\ln
(\frac{\eps(1-\alpha)}{2D_1})\rfloor$, hence the final term is at most
$\eps/2D_1$. The probability that a game completes in any given step is at
least $p$. If we select a time $\tau \geq c/p$ for $c \geq
K+1\geq 1$, then the probability that at most $K$ games are complete is
clearly maximised by taking this probability to be exactly $p$ in all
games. Hence, by Chernoff's bound (see, for example,~\cite[Theorem
2.1]{JLR00}),
\begin{align*}
\e[N(\tau)]&\leq\ {D_2}^K\sum_{k=0}^{K}
\binom{\tau}{k}p^k(1-p)^{\tau-k}+\frac{\eps}{2D_1} \\
&\leq\ e^{K\ln D_2-\frac{(c-K)^2}{2c}}+\frac{\eps}{2D_1} \\
&\leq\ e^{K\ln D_2+K-c/2}+\frac{\eps}{2D_1}.
\end{align*}
Choosing $c=2K\ln(eD_2) +2\ln \frac{2D_1}{\eps}$, we obtain $ \e[N(\tau)] <
\frac{\eps}{D_1} $, where $ \tau =\big\lceil
\frac{3\ln(eD_2)}{p(1-\alpha)}\ln\big(\frac{2D_1}{\eps(1-\alpha)}\big)\big\rceil$.
%\[ \tau\leq \frac{2}{p}\ln\frac{2D_1}{\eps(1-\alpha)}\left(\frac{\ln(eD_2)}{|\ln  \alpha|}+1 \right)+1 < \frac{3\ln(D_2)}{p(1-\alpha)}\ln\frac{(2D_1)}{\eps(1-\alpha)}.\]

We conclude that the gambler's expected return at time $\tau$ is positive.
More importantly, for any initial states $X_0,Y_0\in \Omega$, the expected
distance at time $\tau$ is at most $\eps$ by linearity of expectations, and
so the probability that the chain has not coupled is at most $\eps$. The
mixing time claimed now follows by standard arguments. See, for
example,~\cite{J03}.
\end{proof}

\begin{rem}\label{rem15}
The assumption that the stopping time occurs when the distance changes is
not essential. We clearly cannot dispense with assumption~(ii), or we
cannot bound mixing time. Assumption~(i) may appear a restriction, but
appears to be naturally satisfied in most applications. It seems more
natural than the assumption of bounded stopping time, used in~\cite{HV04}.
Assumption~(i) can easily be replaced by something weaker, for
example by allowing $p$ to vary with time rather than remain constant.
Provided $p\neq 0$ sufficiently often, a similar proof will be valid.
\end{rem}

\begin{rem}\label{rem20}
Let $\gamma=1/(1-\alpha)$. It seems likely that $D_2$ will be small in
comparison to $\gamma$ in most applications, so we might suppose $D_2 <
\gamma < D_1$. The mixing time bound from Theorem~\ref{stopping} can then
be written $O(p^{-1}\gamma\log D_2 \log(D_1/\eps) )$. We may compare this
with the bound which can be derived using~\cite[Corollary~4]{HV04}. This
can be written in similar form as $O(p^{-1}\gamma \log\gamma\log( D_1/\eps)
)$. In such cases we obtain a reduction in the estimate of mixing time by a
factor $\log\gamma/\log D_2$. In the applications below, for example, we
have $D_2=2$ and $\gamma=\Omega(\Delta)$, so the improvement is
$\Omega(\log \Delta)$.
\end{rem}

\begin{rem}
The reason for our improvement on the result of~\cite{HV04} is that the
use of an upper bound on the stopping time, as is done in~\cite{HV04},
will usually underestimate the number of stopping times which occur
in a long interval, and hence the mixing rate.
\end{rem}

\section{Hypergraph independent sets}\label{sec:indsets}

We now use the approach of path coupling via stopping times to prove that
the chain discussed in Section~\ref{sec:intuit} is rapidly mixing. The
metric used in path coupling analyses throughout the paper will be Hamming
distance between the coupled chains. We prove the following theorem.

\begin{thm}\label{indsets}
Let $\lambda,\Delta$ be fixed, and let $\CH$ be a hypergraph such that
$m\geq 2\lambda\Delta+1$. Then the Markov chain $\M(\CH)$ has mixing time
$O(n\log n)$.
\end{thm}
Before commencing the proof itself, we analyse the stopping time
$T$ for this problem.

\subsection{Edge Process}\label{edge}
Let $X_t$ and $Y_t$ be copies of $\M$ which we wish to couple, with
$Y_0=X_0\cup\set{w}$. Let $e$ be any edge containing $w$, with $m=|e|$. We
consider only the times at which some vertex in $e$ is chosen. The progress
of the coupling on $e$ can then be modelled by the following ``game''. We
will call the number of unoccupied vertices in $e$ (excluding $w$)
\emph{units}. At a typical step of the game we have $k$ units, and we
either win the game, win a unit, keep the same state or lose a unit. These
events happen with the following probabilities: we win the game with
probability $1/m$, win a unit with probability at least
$(m-k-1)/(1+\lambda)m$, lose a unit with probability at most $\lambda
k/(1+\lambda)m$ and stay in the same state otherwise. If ever $k=0$, we are
bankrupt and we lose the game. Winning the game models the ``good event''
that the vertex $v$ is chosen and the two chains couple. Losing the game
models the ``bad event'' that the coupling increases the distance to 2. We
wish to know the probability that the game ends in bankruptcy. We are most
interested in the case where $k=1$ initially, which models $e$ being
critical. Note that the value of $k$ in the process on hypergraph
independent sets dominates the value in our model, since we can always
delete (win in the game), but we may not be able to insert (lose in the
game) because the chosen vertex is critical in some other edge.

Let $p_k$ denote the probability that a game is lost, given that
we start with $k$ units. We have the following system of
simultaneous equations.
\begin{align}\label{p-eq1}
    (m-1+2\lambda)p_1 - (m-2)p_2\ &=\ \lambda& \notag\\
   -k\lambda p_{k-1}+(m-k+(k+1)\lambda)p_k-(m-k-1)p_{k+1}\ &=\ 0& (k=2,3,\ldots,m-1)
\end{align}
Adding the equations in~(\ref{p-eq1}) from the $k^\textrm{th}$ onwards
gives
\begin{align}\label{p-eq2}
    (m-1)p_1 + m\lambda p_{m-1}\ &=\ \lambda& \notag\\
   -k\lambda p_{k-1}+(m-k)p_k+ m\lambda p_{m-1}\ &=\ 0& (k=2,3,\ldots,m-1).
\end{align}
Now~(\ref{p-eq2}) is equivalent to~(\ref{p-eq1}), since we have
simply multiplied the coefficient matrix of~(\ref{p-eq1}) by an
upper triangular matrix with all entries 1. This transformation is
clearly nonsingular. We will show by induction that~(\ref{p-eq2}) has solution%
\begin{equation}\label{p-eq3}
    p_k=\frac{\lambda^k-\sum_{i=1}^k \binom{m}{i}p_{m-1}\lambda^{k-i+1}}{\binom{m-1}{k}}
\qquad(k=1,2,\ldots,m-1).
\end{equation}%
When $k=1$, the first equation in~(\ref{p-eq2}) is clearly
satisfied by~(\ref{p-eq3}).  Assume by induction
that~(\ref{p-eq3}) is true for $p_{k-1}$, with $k\geq 2$. Then

\begin{align*}
    p_k\ &=\ \frac{\lambda k}{m-k}\, p_{k-1}-\frac{\lambda m}{m-k}\,p_{m-1}\\
    &=\ \frac{\lambda k}{m-k}\, \frac{\lambda^{k-1}-\sum_{i=1}^{k-1}
    \binom{m}{i}\lambda^{k-i}p_{m-1}}{\binom{m-1}{k-1}}-\frac{\lambda m}{m-k}\,p_{m-1}\\
    &=\ \frac{\lambda^k-\sum_{i=1}^{k-1} \binom{m}{i}\lambda^{k-i+1}p_{m-1}}{\binom{m-1}{k}}
    - \frac{\binom{m}{k}}{\binom{m-1}{k}}\,\lambda p_{m-1}\\
    &=\ \frac{\lambda^k-\sum_{i=1}^{k} \binom{m}{i}\lambda^{k-i+1}p_{m-1}}{\binom{m-1}{k}},
\end{align*}
continuing the induction. For consistency, we must clearly have
\begin{align}
    p_{m-1}\ &=\ \frac{\lambda^{m-1}-\sum_{i=1}^{m-1}
    \binom{m}{i}\lambda^{m-i}p_{m-1}}{\binom{m-1}{m-1}}
    \ =\ \lambda^{m-1}-\big((1+\lambda)^m-1-\lambda^m\big)p_{m-1},\notag\\
    \textrm{i.e.}\quad p_{m-1}\ &=\ \frac{\lambda^{m-1}}{(1+\lambda)^m-\lambda^m}.\label{p-eq4}
\end{align}%
Using (\ref{p-eq4}),~(\ref{p-eq3}) can be rewritten
\begin{equation}\label{p-eq5}
    p_k\ =\ \frac{1}{\binom{m-1}{k}}\bigg(\lambda^k-\frac{\sum_{i=1}^k
    \binom{m}{i}\lambda^{m+k-i}}{(1+\lambda)^m-\lambda^m}\bigg)
    = \frac{\sum_{i=k+1}^m \binom{m}{i}\lambda^{m+k-i}}
    {\big((1+\lambda)^m-\lambda^m\big)\binom{m-1}{k}}
\qquad(k=1,2,\ldots,m-1).
\end{equation}%
In particular
\begin{equation}\label{p-eq6}
p_1= \frac{\lambda}{m-1}\left(1-\frac{m\lambda^{m-1}}
{(1+\lambda)^m-\lambda^m}\right).
\end{equation}

\subsection{The expected distance between $X_T$ and $Y_T$}
The stopping time for the pair of chains $X_t$ and $Y_t$ will be when the
distance between them changes, in other words either a good or bad event
occurs. The probability that we observe the bad event on a particular edge
$e$ with $w\in e$ is at most $p_k$ as calculated above. Let $\xi_t$ denote
the number of empty vertices in $e$ at time $t$ when the process is started
with $\xi_0=k$. Now $\xi_t$ can never reach 0 without first reaching $k-1$
and, since the process is Markovian, it follows that
\[p_k= \Pr(\exists t\, \xi_t=0 | \xi_0=k)
= \Pr(\exists t\,\xi_t=0 | \xi_{s}=k-1)\Pr(\exists s\,\xi_{s}=k-1 | \xi_{0}=k)
< p_{k-1}.\]%
Since $w$ is in at most $\Delta$ edges, the probability that we observe the
bad event on any edge is at most $\Delta p_1$. The probability that the
stopping time ends with the good event is therefore at least $1-\Delta
p_1$.  The path coupling
calculation is then%
\[\e[\mathrm{d}(X_T,Y_T)]\leq 2\Delta p_1.\]
This is required to be less than 1 in ordered to apply
Theorem~\ref{stopping}. If $m \geq 2\lambda\Delta+1$, then by~(\ref{p-eq6})
\[2\Delta p_1=1-\frac{(2\lambda\Delta+1)\lambda^{2\lambda\Delta}}
{(1+\lambda)^{2\lambda\Delta+1}-\lambda^{2\lambda\Delta+1}}.\]

\begin{proof}[Proof of Theorem~\ref{indsets}]
The above work puts us in a position to apply
Theorem~\ref{stopping}. Let $m\geq 2\lambda\Delta+1$. Then for
$\M(\CH)$ we have
\begin{enumerate}
\item $\Pr(\mathrm{d}(X_t,Y_t)\neq1| \mathrm{d}(X_{t-1},Y_{t-1})=1)\geq
\frac{1}{n}$ for all $t$, and \item $\e[\mathrm{d}(X_T,Y_T)]\ <\ 1\, -\,
\dfrac{(2\lambda\Delta+1)\lambda^{2\lambda\Delta}}
{(1+\lambda)^{2\lambda\Delta+1}-\lambda^{2\lambda\Delta+1}}.$
\end{enumerate}
Also for $\M(\CH)$ we have $D_1=n$ and $D_2=2$. Hence by
Theorem~\ref{stopping}, $\M(\CH)$ mixes in time
 \[\tau(\eps)\leq
6n \frac{(1+\lambda)^{2\lambda\Delta+1}-\lambda^{2\lambda\Delta+1}}
{(2\lambda\Delta+1)\lambda^{2\lambda\Delta}}\ln\Big(n\eps^{-1}\frac{(1+\lambda)^{2\lambda\Delta+1}-\lambda^{2\lambda\Delta+1}}
{(2\lambda\Delta+1)\lambda^{2\lambda\Delta}}\Big).\]

This is $O(n\log n)$ for fixed $\lambda,\Delta$.
\end{proof}
\begin{rem}
In the most important case, $\lambda=1$, we require $m\geq
2\Delta+1$. This does not include the case $m=3$, $\Delta=3$
considered in~\cite{DG00a}. We have attempted to improve the bound by employing the chain proposed by
Dyer and Greenhill in~\cite[Section~4]{DG00a}. However, this gives only a
marginal improvement. For large $\lambda\Delta$, we obtain convergence for
$m \geq 2\lambda\Delta+\tfrac{1}{2}+o(1)$. For $\lambda=1$, this gives a
better bound on mixing time for $m=2\Delta+1$, with dependence on $\Delta$
similar to Remark~\ref{rem05} below, but does not even achieve mixing for
$m=2\Delta$.
We omit the details in
order to deal with the Glauber dynamics, and to simplify the
analysis.
\end{rem}
\begin{rem}\label{rem05}
The terms in the running time which are exponential in
$\lambda,\Delta$ would disappear if we instead took graphs for which
$m\geq 2\lambda\Delta+2$. In this case the running time would be
\[\tau(\eps)\leq
 6(2\lambda\Delta+1) n \ln(n\eps^{-1}(2\lambda\Delta+1))\leq  12(2\lambda\Delta+1) n \ln(n\eps^{-1}).\]
Furthermore, if we took graphs such that
$m>(2+\delta)\lambda\Delta$, for some $\delta>0$, then the running
time would no longer depend on $\lambda,\Delta$ at all, but would be
$\tau(\eps)\leq c_\delta n \ln (n\eps^{-1})$ for some constant
$c_\delta$.
\end{rem}
\begin{rem}\label{rem10}
It seems that path coupling cannot show anything better than $m$ linear in
$\lambda\Delta$. Suppose the initial configuration has edges
$\{w,v_1,\ldots,v_{m-2},x_i\}$ for $i=1,\ldots,\Delta$, with
$w,v_1,\ldots,v_{m-2}\in X_0$, $x_1,\ldots,x_\Delta\not\in X_0$ and $w$ the
change vertex. Consider the first step where any vertex changes state. Let
$\mu=(1+\lambda)(m-1+\Delta)$. The good event occurs with probability
$(1+\lambda)/\mu$, insertion of a critical vertex with probability
$\lambda\Delta/\mu$, and deletion of a non-critical vertex with probability
$(m-1)/\mu$. We therefore need $(m-1)+(1+\lambda)\geq\lambda\Delta$, i.e.
$m\geq \lambda(\Delta-1)$, to show convergence by path coupling.
\end{rem}
\begin{rem}
It seems we could improve our bound $m\geq 2\lambda\Delta+1$ for rapid
mixing of the Glauber dynamics somewhat if we could analyse the process on
all edges simultaneously. Examination of the extreme cases, where all edges
adjacent to $w$ are otherwise independent, or where they are dependent
except for one vertex (as in Remark~\ref{rem10}), suggests that improvement
to $(1+o(1))\lambda\Delta$ may be possible, where the $o(1)$ is relative to
$\lambda\Delta$. However, the analysis in the general case seems difficult,
since edges can intersect arbitrarily.
\end{rem}

\section{Hardness results for independent sets}\label{sec:hard}
We have established that  the number of independent sets of a hypergraph
can be approximated efficiently using the Markov Chain Monte Carlo
technique for hypergraphs with edge size linear in $\Delta$. We show next
that exact counting is unlikely to be possible, and that our approximation
scheme cannot be extended to cover all hypergraphs with edge size
$\Omega(\log\Delta)$.

\subsection{\#P-completeness}\label{sec:hardnump}
We show that the exact counting problem is \#P-Complete except in trivial
cases.
\begin{thm}\label{h-thm10}
Let $\CG(m,\Delta)$ be the class of uniform hypergraphs with minimum edge
size $m\geq 3$ and maximum degree $\Delta$. Computing the number of independent
sets of hypergraphs in $\CG(m,\Delta)$ is \#P-complete if $\Delta \geq 2$.
If $\Delta \leq 1$, it is in P.
\end{thm}
\begin{proof}
Since $m$ is the minimum edge size, we will assume $m\geq 3$. The cases
$\Delta=0,1$ are trivially in P. As discussed in Section~\ref{sec:intuit},
independent sets in a hypergraph with $\Delta=2$ correspond to edge covers
in a graph. Counting these is \#P-complete, even for graphs with
arbitrarily large minimum degree. This is stated in~\cite{BD97} but without
proof, so we provide a proof in Appendix~\ref{app1}. We now consider
$\Delta \geq 3$. (The case $m=\Delta=3$ is discussed in~\cite{DG00a}.) Take
a graph $G=(V,E)$, and construct a hypergraph $\CG=(\CV,\CE)$ by
``extending'' each edge $e=\set{v_1,v_2}\in E$ to an edge
$e^+=\set{v_1,u^e_1,\ldots, u^e_{m-2},v_2}\in \CE$. Observe that, for each
independent set $I$ of $G$ and edge $e\in E$, there are $2^{m-2}-1$
independent assignments to $u^e_1,\ldots, u^e_{m-2}$ if $v_1,v_2\in I$ and
$2^{m-2}$ otherwise. This is equivalent to evaluating the partition
function of a \emph{weighted $H$-colouring} problem~\cite{BG04,DG00b} on
$G$, with weight matrix
\[ A=\begin{bmatrix} 2^{m-2} & 2^{m-2}\\ 2^{m-2} & 2^{m-2}-1 \end{bmatrix}.\]
The \#P-completeness of $H$-colouring with this weight matrix follows
either directly from \cite{BG04} or indirectly from \cite[Corollary
3.2]{DG00b}. The degree bound $\Delta=3$ follows from \cite[Theorem
5.1]{DG00b}, on noting that $A$ is nonsingular.
\end{proof}
\subsection{Approximation hardness}
We now show that unless NP\/=\/RP, there can be no \emph{fpras}
for the number of independent sets of all hypergraphs with edge size
$\Omega(\log\Delta)$.

Let $G=(V,E)$, with $|V|=n$, be a graph with maximum degree
$\Delta$ and $N_i$ independent sets of size $i$
($i=0,2,\ldots,n$). For $\lambda>0$ let $Z_G(\lambda)=\sum_{i=0}^n
N_i \lambda^i$ define the \emph{hard core partition function}. The
following is a combination of results in Luby and
Vigoda~\cite{LV99} and Berman and Karpinski~\cite{BK03}.
\begin{thm}\label{h-thm20}
If $\lambda > 694/\Delta$, there is no \emph{fpras} for
$Z_G(\lambda)$ unless NP\/=\/RP.
\end{thm}
\begin{proof}
Let $\eps$ be a constant such that the size of the largest
independent set in a graph of maximum degree 4 cannot be
approximated to within a ratio $(1+\eps)$ unless P\,=\,NP. Berman
and Karpinski~\cite{BK03} show that $\eps\geq 1/49$. Luby and
Vigoda~\cite[Theorem~4]{LV99} prove the hardness of approximating
$Z_G(\lambda)$ if $\lambda > c/\Delta$ for any $c > 20\ln
2\,(1+\eps)/\eps$.\footnote{The expression in~\cite{LV99} omits
the $\ln 2$ term} Together, these two results give the theorem.
\end{proof}
We note that Theorem~\ref{h-thm20} could probably be strengthened using the
approach of~\cite{DFJ99}. However, this has yet to be done.
\begin{thm}\label{h-thm30}
Unless NP\/=\/RP, there is no \emph{fpras} for counting
independent sets in hypergraphs with maximum degree $\Delta$ and
minimum edge size $m< 2\lg(1+\Delta/694)-1=\Omega(\log\Delta)$.
\end{thm}
\begin{proof}
Given a graph $G=(V,E)$ with maximum degree $\Delta$, we construct
a hypergraph $\CH=(\CV,\CE)$ as follows. Let $k=\lceil m/2\rceil$.
For each $v\in V$, let $W_v=\set{w_{v1},w_{v2},\ldots,w_{vk}}$ and
$\CV=\bigcup_{v\in V} W_v$. For each edge $e=\set{u,v}\in E$, let
$S_e=W_u\cup W_v$, and let $\CE=\set{S_e:e\in E}$. It is clear
that $\CH$ has maximum vertex degree $\Delta$ and every edge has
size $2k\geq m$.

An independent set $\CI$ in $\CH$ corresponds to a unique
independent set $I$ in $G$ as follows. If $S_v\subseteq \CI$, then
$v\in I$, otherwise $v\notin I$. Clearly $\CI$ independent in
$\CH$ implies $I$ independent in $G$. Note that for each $v\notin
I$, there are $(2^k-1)$ possible subsets of $S_v$ which may be in
$\CI$. Thus, if $\CN$ is the number of independent sets in $\CH$,
\begin{equation*}%\label{h-eq1}
    \CN=\sum_{i=0}^n N_i(2^k-1)^{n-i}= (2^k-1)^n \sum_{i=0}^n N_i(2^k-1)^{-i}
    =(2^k-1)^n Z_G(1/(2^k-1)).
\end{equation*}
Thus approximating $\CN$ is equivalent to approximating
$Z_G(\lambda)$ with $\lambda=1/(2^k-1)$. But, by
Theorem~\ref{h-thm20}, this will be hard if $1/(2^k-1)>694/\Delta$.
This gives $k<\lg(1+\Delta/694)$, which holds whenever
$m<2\lg(1+\Delta/694)-1$.
\end{proof}

% ----------------------------------------------------------------
\section{Hypergraph colouring}\label{sec:colour}
We now consider Glauber dynamics on the set of proper colourings of a
hypergraph. Again our hypergraph $\CH$ will have maximum degree $\Delta$,
minimum edge size $m$, and we will have a set of $q$ colours. A colouring
of the vertices of $\CH$ is proper if no edge is monochromatic. Let
$\Omega'(\CH)$ be the set of all proper $q$-colourings of $\CH$. We define
the Markov chain $\C(\CH)$ with state space $\Omega'(\CH)$ by the following
transition process. If the state of $\C$ at time $t$ is $X_t$, the state at
$t+1$ is determined by
\begin{enumerate}
\item selecting a vertex $v\in \CV$ and a colour $k\in\{1,2,\ldots,q\}$
uniformly at random, \item let $X'_t$ be the colouring obtained by
recolouring $v$
colour $k$ \item if $X'_t$ is a proper colouring let $X_{t+1}=X'_t$\\
otherwise let $X_{t+1}=X_t$.
\end{enumerate}
This chain is easily shown to be ergodic with the uniform
stationary distribution. Again we will use Theorem~\ref{stopping}
to prove rapid mixing of this chain under certain conditions,
however first we will examine the chain using standard path
coupling techniques.

\begin{thm}\label{colsmge4}
For $m\geq 4$, $q> \Delta$, the Markov chain $\C(\CH)$ mixes in time
$O(n\log n)$.
\end{thm}
\begin{proof}
Suppose that two copies of $\C(\CH)$, $X_0$ and $Y_0$ say, start at
distance one apart, i.e. they differ in only one vertex $w$. Suppose that
the number of colours available for recolouring $w$ is $q-k$, then the
probability of the two copies of the chain coupling in one step is
$\frac{q-k}{nq}$. The distance between the two chains can only increase (to
2) if we select a vertex $v$ and recolour it with a colour that is
permitted in one copy of the chain only. For this to happen, there must be
an edge containing $v$ and $w$ such that the other vertices in this edge
are all either red and we have chosen red for $v$, or blue and we have
chosen blue for $v$. Hence there can be at most one vertex on each edge,
and one colour for that vertex, such that the chains diverge if we select
that vertex and colour. Furthermore, for each of the $k$ unavailable
colours there must be an edge containing $w$ which, apart from $w$ itself,
is monochromatic in the forbidden colour, so on these edges there are no
vertices whose selection can cause the chains to diverge. Hence the
probability that the distance increases to 2 in one step is at most
$\frac{\Delta-k}{nq}$. The path coupling calculation is therefore
\[ \e[\mathrm{d}(X_1,Y_1)]\leq 1-\frac{q-k}{nq}+\frac{\Delta-k}{nq}.\]
If $q\geq\Delta+1$ then $\e[\mathrm{d}(X_1,Y_1)]\leq 1-1/nq$, and therefore
by the path coupling theorem the mixing time is
\[ \tau(\eps)\leq nq \ln (n \eps^{-1}).\vspace{-\baselineskip}\]
\end{proof}

This analysis leaves little room for improvement in the case $m\geq 4$,
indeed it is not clear whether the Markov chain described is even ergodic
for $q\leq \Delta$. The following simple construction does show that the
chain is not in general ergodic if $q\leq \frac{\Delta}{m}+1$. Let
$q=\frac{\Delta }{m }+1$, and take a hypergraph $\CH$ on $q(m-1)$ vertices.
We will group the vertices into $q$ groups $\CV=\CV_1,\CV_2,\ldots,\CV_q$,
each of size $m-1$. Then the edge set of $\CH$ is $E=\{\{v\}\cup \CV_j:v\in
\CV, v\not\in \CV_j\}$. The degree of each vertex is
$(q-1)+(q-1)(m-1)=\Delta$. If we now colour each group $\CV_j$ a different
colour, we obtain $q!$ distinct colourings, but for each of these the
Markov chain is frozen (no transition is valid).

The case $m=2$ is graph colouring and has been extensively studied. See,
for example,~\cite{DFHV04}. This leaves the case $m=3$, hypergraphs with 3
vertices in each edge. The standard path coupling argument, as in
Theorem~\ref{colsmge4}, only shows rapid mixing for $q\geq 2\Delta$, since
there may be two vertices in each edge that can be selected and lead to a
divergence of the two chains. This occurs if, of the two vertices in an
edge which are not $w$, one is coloured red and the other blue. However, we
can do better using Theorem~\ref{stopping}. We will need the following
technical Lemma.
\begin{lem}\label{dws}
Let $ \varphi(d)=1- d(1-e^{-(q-\Delta+d)t/Mq})/(q-\Delta+d)$. For
all $t\geq 0$ and all $d\geq 1$, $\varphi(d)\geq\varphi(1)^d$.
\end{lem}
\begin{proof}
Let $\kappa=q-\Delta> 0$,
$x=t/Mq\geq 0$.
We wish to show that\vspace{-1.5ex}%
\[ \psi(x)= \varphi(d)-\varphi(1)^d=
1- d(1-e^{-(\kappa+d)x})/(\kappa+d) - \big(1-
(1-e^{-(\kappa+1)x})/(\kappa+1)\big)^d \geq 0.\]%
Since $\psi(0)=0$, it suffices to show that $\psi(x)$ is
increasing for all $x \geq 0$. But
\begin{align*}
    \psi'(x)\ &=\ -de^{-(\kappa+d)x}+de^{-(\kappa+1)x}\big(1-
    (1-e^{-(\kappa+1)x})/(\kappa+1)\big)^{d-1}\\[1ex]
    &=\ de^{-(\kappa+1)x}\big(\big(1-(1-e^{-(\kappa+1)x})/(\kappa+1)\big)^{d-1}
    -e^{-(d-1)x}\big),
\end{align*}
so it suffices to show $1-(1-e^{-(\kappa+1)x})/(\kappa+1) \geq
e^{-x}$. Let $\zeta(x)=1-e^{-x}-(1-e^{-(\kappa+1)x})/(\kappa+1)$.
Then $\zeta(0)=0$, so we need only show that $\zeta(x)$ is
increasing for all $x\in(0,\infty)$. But
\begin{equation*}
    \zeta'(x)\ =\ e^{-x}-e^{-(\kappa+1)x}
    =\ e^{-x}(1-e^{-\kappa x})\ \geq\ 0,
\end{equation*}
for $x \geq 0$.
\end{proof}\vspace{0ex}

\begin{thm}\label{colouring}
There exists $\Delta_0$ such that, if $\CH$ is a 3-uniform hypergraph with
maximum degree $\Delta>\Delta_0$ and $q \geq 1.65 \Delta$, the Markov chain
$\C(\CH)$ mixes rapidly.
\end{thm}

\begin{proof}
We choose $\Delta_0$ large enough that all the approximations below are
valid. We couple two copies of this chain using the identity coupling. Let
$X$ and $Y$ be two copies of $\C(\CH)$ such that $X_0$ and $Y_0$ differ
only at a single vertex $w$. As before, we will examine the stopping time
$T$ at which $\mathrm{d}(X_T,Y_T)\neq 1$ for the first time, and show that
$\e[\mathrm{d}(X_T,Y_T)]< 1$. We assume that $w$ is coloured blue in $X_0$
and red in $Y_0$. We will call any other colour \emph{neutral}. Let
$\Gamma(w)$ denote the set of vertices of $\CH$ that share an edge with
$w$. We will only consider transitions in which either $w$ or a vertex in
$\Gamma(w)$ is selected, since any transition which involves any other
vertex will not change the distance between $X$ and $Y$. Let
$M=|\Gamma(w)|+1$. We will first assume that none of the edges containing
$w$ is otherwise monochromatic, and hence that all colours are available
for recolouring $w$. We will deal with other cases later. Let $S_t$ denote
the event that $T=t$ and $\mathrm{d}(X_T,Y_T)=0$, which we will call
\emph{success}. The bad event we will call \emph{failure}.

The probability that the two chains couple in any one step is
$\frac{q}{Mq}$. For each $v\in \Gamma(w)$, let $\beta_{v,t}$ be an
indicator variable which takes value 1 if $v$ is either red or blue after
$t$ steps of the chain, and takes value 0 otherwise. We describe a choice
of vertex $v\in\Gamma(w)$ and colour $c\in \{\textrm{red,blue}\}$ at step
$t$ as `bad' if there is an edge containing $v$ and $w$ whose other vertex
is currently coloured $c$, and let $B_t$ denote the number of bad choices
at time $t$. The probability of failure in step $t$ is therefore
$\frac{B_t}{Mq}$. For each $v\in \Gamma(w)$ let $d_v$ be the number of
edges which contain both $v$ and $w$. Then $B_t\leq \sum_{v\in\Gamma(w)}
d_v\beta_{v,t}$.  Now, using $\approx$ to imply equality up to a factor
$1+o_\Delta(1)$,

\begin{align}
 \Pr(S_t)&=\ \e\bigg[\prod_{j=0}^{t-1}\left(1-\frac{1}{M}-\frac{B_t}{Mq}\right)
 \frac{1}{M}\bigg]\ \
 \approx\  \ \frac{1}{M}\e\Big[e^{-\sum_{j=0}^{t-1}(\frac{1}{M}+\frac{B_t}{Mq})}\Big]\notag\\[1ex]
& \geq\ \frac{1}{M}\e\Big[e^{-\sum_{j=0}^{t-1}(\frac{1}{M}+\frac{\sum_{w\in\Gamma(v)}
 d_w\beta_{w,t}}{Mq})}\Big]\
% &=\ \frac{e^{-t/M}}{M}e^{-\sum_{w\in\Gamma(v)} \frac{ d_w}{Mq}\sum_{j=0}^{t-1}\beta_{w,t}}\\
 =\ \frac{e^{-t/M}}{M}\e\Big[\prod_{w\in\Gamma(v)}e^{- \frac{
d_w}{Mq}\sum_{j=0}^{t-1}\beta_{w,t}}\Big].\label{prs}
\end{align}

We will now study the properties of $\beta_{v,t}$, with a view to analysing
$\e[e^{- \frac{ d_v}{Mq}\sum_{j=0}^{t-1}\beta_{v,t}}]$. Note that the
probability has not coupled or diverged by time $40\Delta$
is at most%
\[(1-1/M)^{20\Delta}\leq (1-1/2\Delta)^{40\Delta} \leq e^{-20}< 10^{-8},\]
so we consider times only up to $40\Delta$. Let $t_v = \max \{t:t< 40\Delta
\textrm{ and }\beta_{v,t}=1\}$. If $v$ starts out either red or blue, the
probability that it is recoloured to a neutral colour in each step is at
least $(q-\Delta-2)/Mq$. Also, the probability that it becomes red or blue
before time $40\Delta$ is at most $80\Delta/Mq$. Hence
\begin{align*}
\Pr(t_w > t ) &\leq\ \left( 1- \frac{q-\Delta-2}{Mq}\right)^t + \frac {80\Delta}{Mq}\\
% &\approx\ \left( 1- \frac{q-\Delta}{Mq}\right)^t\\
& \approx\ e^{-\frac{q-\Delta}{Mq}t},
\end{align*}
since the second term is $O(1/\Delta)$ and small compared to the first,
which is $\Omega(1)$ for $t\leq 40\Delta$. Now we can bound
$\sum_{j=0}^{t-1}\beta_{v,t}$ by the minimum of $t$ and $t_v$, an
exponentially distributed random variable with parameter
$\frac{q-\Delta}{Mq}$. We are in a position to bound $\e[e^{- \frac{
d_v}{Mq}\sum_{j=0}^{t-1}\beta_{v,t}}]$ as follows.
\begin{align*}
\e[e^{- \frac{ d_v}{Mq}\sum_{j=0}^{t-1}\beta_{v,t}}] & \geq
\ \sum_{j=0}^{t} \Pr(t_v=j)e^{-\frac{ d_v}{Mq}j} + \Pr(t_v > t )
e^{-\frac{ d_v}{Mq}t}\\
&\approx \ \int_{0}^{t} \frac{q-\Delta}{Mq}e^{- \frac{q-\Delta}{Mq}x}
e^{-\frac{d_v}{Mq} x} dx + e^{- \frac{q-\Delta+d_v}{Mq}t}\\
&=\ 1- \frac{d_v}{q-\Delta+d_v}\left(1-e^{-\frac{q-\Delta+d_v}{Mq}
t}\right).
\end{align*}
Inserting this into Equation~(\ref{prs}), we get
\[\e[\Pr(S_t)]\ \geq\ \frac{e^{-t/M}}{M}\prod_{v\in\Gamma(w)}
\left(1- \frac{d_v}{q-\Delta+d_v}\left(1-e^{-\frac{q-\Delta+d_v}{Mq}
t}\right)\right)=\frac{e^{-t/M}}{M}\prod_{v\in \Gamma(w)} \varphi(d_v),\]%
where $ \varphi(d)$ was defined in Lemma~\ref{dws}. Since $\Sigma_{v\in \Gamma(w)}d_v=2\Delta$, Lemma~\ref{dws} implies that for all $t\geq 0$,
\[ \prod_{v\in \Gamma(w)} \varphi(d_v) \geq
\varphi(1)^{2\Delta}.\]
Hence, for $t \leq 40\Delta$,
\begin{align*}
\e[\Pr(S_t)]\ &\geq\ \frac{e^{-t/M}}{M} \left(1- \frac{1}{q-\Delta+1}\left(1-e^{-\frac{q-\Delta+1}{Mq} t}\right)\right)^{2\Delta}\\
&\approx \ \frac{1}{M}
e^{-\frac{t}{M}-\frac{2\Delta}{q-\Delta}(1-e^{-\frac{q-\Delta}{Mq}t})}.
\end{align*}
Finally, noting that $\Pr(\mathrm{d}(X_T,Y_T)=0)=\sum_{t=0}^{\infty}
\Pr(S_t)$ by linearity of expectation, we have
\begin{align*}
\Pr(\mathrm{d}(X_T,Y_T)=0)\ &\geq\ \int_0^{40\Delta} \frac{1}{M} e^{-\frac{t}{M}-\frac{2\Delta}{q-\Delta}(1-e^{-\frac{q-\Delta}{Mq}t})}dt\\
&=\ \int_0^{\frac{40 \Delta}{M}}
e^{-z-\frac{2\Delta}{q-\Delta}(1-e^{-\frac{q-\Delta}{q}z})}dz
\end{align*}
If we now substitute $q=1.65 \Delta$ and $M\leq 2\Delta$, we see
that
\[\Pr(\mathrm{d}(X_T,Y_T)=0)\ \geq\ \int_0^{20} e^{-z-3.077(1-e^{-0.3941z})}dz\ >\
0.5003.\] Since $\mathrm{d}(X_T,Y_T)\in \{1,2\}$, it follows that
$\e[\mathrm{d}(X_T,Y_T)]< 0.9994$ and we can apply Theorem~\ref{stopping}.
This yields the claimed result.

We have assumed that all colours are available for recolouring $w$ at every
step. This will not be the case if there is any edge $e$ adjacent to $w$
for which $e\setminus\{w\}$ is monochromatic. Let us call such an edge
\emph{blocking}, and suppose there are $\rho_t$ blocking edges at time $t$.
Note that the failure cannot occur on a blocking edge. The total number
$\rho'$ of blocking edges created during time $40\Delta$ is at
most $\ln\Delta$ since%
\[ \Pr(\rho' \geq \ln\Delta)\ \leq\
\binom{40\Delta}{\ln\Delta}\Big(\frac{1}{q}\Big)^{\ln\Delta}\ \leq\
\Big(\frac{40}{\ln\Delta}\Big)^{\ln\Delta}\ =\ O\Big(\frac{1}{\Delta^\gamma}\Big)\]%
for every constant $\gamma>0$. Since $\ln\Delta$ is negligible in
comparison with $q$ and $M$, these do not affect the probability estimates
in the proof above. Thus we may assume that all blocking edges exist
initially. We may further assume these persist until termination, so
$\rho_t=\rho_0=\rho$ for all $t$. This can only decrease the probability of
success. We now observe that this is no worse in our analysis than taking
$\Delta'=\Delta-\rho$ and $q'=q-\rho$. Let us formally define $M'=Mq/q'$.
Then the conditional success probability is $(q-\rho)/Mq=1/M'$, and the
recolouring probability at each step is at least
$(q-\Delta)/Mq=(q'-\Delta')/M'q'$. The analysis now proceeds as before.
Since $M'$ plays no part in the final condition, we finally require $q'
\geq 1.65\Delta'$, i.e. $q \geq 1.65\Delta-0.65\rho$. This is clearly a
weaker condition than $q \geq 1.65\Delta$.
\end{proof}
\begin{rem}
If we let $\beta=(q-\Delta)/q$ then, as $\Delta_0\rightarrow\infty$, the
analysis can be tightened slightly to work for $\beta>\beta^*$, where
$\beta^*$ is the root of the
equation%
\[ \int_0^{\infty}
e^{-z-\frac{2(1-\beta)}{\beta}(1-e^{-\beta z})}dz\ =\ \tfrac{1}{2}.\]%
The integral can be expanded, by parts integration, as an infinite series
to give an alternative equation%
\[ \sum_{i=0}^\infty \frac{(-2)^i(1-\beta)^i}{\prod_{j=0}^i(1+j\beta)}\
=\ \tfrac{1}{2}.\]%
This has root $\beta^*=0.392729$, giving $q > 1.64671$.
\end{rem}
\begin{rem}
A route to improving our bound on $q$ would be to consider the changes in
the numbers of colours available at each vertex of $\Gamma(w)$ during the
process. We make the pessimistic assumption that this is always $q-\Delta$
but, while this could be true initially, we would expect more colours to
become available later on. A proper analysis of this effect seems more
difficult, however, because $\Theta(\Delta^2)$ vertices are now involved,
and the edges containing them may intersect.
\end{rem}
\section{Hardness results for colouring}\label{sec:hardcol}
\subsection{\#P-completeness}
Again we show that exact counting is \#P-complete except in the few cases
where it is clearly in P. Let $\CG(m,\Delta)$ be as in
Theorem~\ref{h-thm10}.
\begin{thm}\label{h-thm40}
Computing the number of $q$-colourings of hypergraphs in $\CG(m,\Delta)$ is
\#P-complete if $\Delta,q > 1$. If $\Delta \leq 1$ or $q\leq 1$ it is in P.
\end{thm}
\begin{proof}
Again we assume $m\geq 3$. The cases $\Delta\leq 1$, $q\leq 1$ are
trivially in P. The case $\Delta=2$ corresponds to counting edge
$q$-colourings of graphs in which no vertex is monochromatic. We call an
edge colouring with no monochromatic vertex a \emph{weak edge colouring}.
Counting weak edge colourings is \#P-complete for graphs of arbitrarily
large minimum degree. We give a proof in Appendix~\ref{app2}.

For $\Delta \geq 3$, $q\geq 2$, we use the construction from the proof of
Theorem~\ref{h-thm10}. For each colouring $X:V\rightarrow \{1,2,\ldots,q\}$
of $G$ and edge $e\in E$, there are $q^{m-2}-1$ permitted colourings of
$u^e_1,\ldots, u^e_{m-2}$ if $X(v_1)=X(v_2)$ and $q^{m-2}$ otherwise. The
corresponding $H$-colouring problem has the following $q\times q$ weight
matrix:
\[ A=\begin{bmatrix} q^{m-2}-1 & q^{m-2} & \cdots & q^{m-2}\\
q^{m-2} & q^{m-2}-1& \cdots& q^{m-2}\\
\vdots&\ \ \ddots&&\vdots\\
q^{m-2} & q^{m-2}& \cdots & q^{m-2}-1\end{bmatrix}.\]
The \#P-completeness of $H$-colouring with this weight matrix, and
the bound $\Delta=3$, follow as in Theorem~\ref{h-thm10},
since $A$ is again nonsingular.
\end{proof}
\subsection{Hardness of Approximation}
Again let $\CG(m,\Delta)$ be as defined in Theorem~\ref{h-thm10}. Our
result, Corollary~\ref{h-cor10}, follows directly from the following
NP-completeness proof.
\begin{thm}\label{h-thm50}
Determining whether a hypergraph in $\CG(m,\Delta)$ has any $q$-colouring
is NP-complete for any $m > 1$ and $2 < q \leq(1-1/m)\Delta^{1/(m-1)}$.
\end{thm}
\begin{proof}

If $m=2$, this is graph colouring, and the result follows
from~\cite[Theorem 1.4]{EHK98}. (See also~\cite{MR01}.) For $m\geq 3$, we
use the following reduction from graph colouring. Let $G=(V,E)$ be a graph
with degree $\Delta_G$, and $2 < q \leq 3\Delta_G/4$. Without loss, we may assume
$\Delta_G= \lceil 4q/3\rceil$. Colouring $G$ with $q$ colours is
NP-complete~\cite{EHK98}. For each edge $e=\set{v_1,v_2}\in E$, let
$S_i^e=\set{u^e_{i1},u^e_{i2},\ldots,u^e_{im}}$ $(i=1,2,\ldots,q)$ and
$\CV_0^e=\bigcup_{i=1}^q S_i^e$. Let $\CE_0^e$ comprise all subsets of
$\CV_0^e$ of size $m$ other than $S_i^e$ $(i=1,2,\ldots,q)$. We claim that
any proper $q$-colouring of the hypergraph $\CH_0^e=(\CV_0^e,\CE_0^e)$ must
assign the same colour to all $u^e_{ij}\in S_i^e$ $(j=1,2,\ldots,m)$, and a
different colour for each $i=1,\ldots,q$. The claim holds since there must
be some colour class of size at least $m$, since there are $q$ colours and
$mq$ vertices. If there was a colour class of size greater than $m$, at
least one of its subsets of size $m$ would be a monochromatic edge. Thus
there must be exactly $q$ colour classes, each of size $m$. If these are
not the $S_i^e$ $(i=1,2,\ldots,q)$, again there is a monochromatic subset
of size $m$ which is an edge. Clearly, by symmetry, any assignment of the
$q$ colours to the $q$ classes $S_i^e$ is permissible.

Let $\CV^e=\CV_0^e\cup\set{v_1,v_2}$, and add the edges
$\set{v_1,u^e_{i2},\ldots,u^e_{im}}$ $(i=1,\ldots,\lfloor q/2\rfloor)$ and
$\set{v_2,u^e_{i2},\ldots,u^e_{im}}$ $(i=\lfloor q/2\rfloor+1,\ldots,q)$ to
$\CE_0^e$ to give $\CE^e$. We claim that, in any proper $q$-colouring of
the hypergraph $\CH^e=(\CV^e,\CE^e)$, $v_1$ and $v_2$ must receive
different colours. The claim holds since $v_1$ can have any colour
different from all $S_i^e$ $(i=1,2,\ldots,\lfloor q/2\rfloor)$, and $v_2$
any colour different from all $S_i^e$ $(i=\lfloor q/2\rfloor+1,\ldots,q)$.
But these permitted colour sets for $v_1$ and $v_2$ are disjoint. Also,
given any colours for $v_1$ and $v_2$, there are $\lfloor q/2\rfloor\lceil
q/2\rceil(q-2)!>0$ colourings of $\CH^e$. Thus we may use $\CH^e$ to
simulate the edge $e\in E$. Thus we set $\CV= \bigcup_{e\in E}\CV^e$,
$\CE=\bigcup_{e\in E}\CE^e$ and consider the hypergraph $\CH=(\CV,\CE)$.
Then $\CH$ is $q$-colourable if and only if $G$ is $q$-colourable.

The maximum degree in $\CH$ of any $u_{ij}$ is $\binom{mq}{m-1} \geq
eq^{m-1}$. The degree in $\CH$ of each $v\in V$ is at most $\Delta_G\lceil
q/2\rceil \leq (4q+2)(q+1)/6 < 2q^2$. Thus $\Delta=\binom{mq}{m-1}\geq
(mq/(m-1))^{m-1}$, and hence $q\leq (1-1/m)\Delta^{1/(m-1)}$.
\end{proof}
\begin{cor}\label{h-cor10}
Unless NP\/=\/RP, there is no \emph{fpras} for counting $q$-colourings of a
hypergraphs with maximum degree $\Delta$ and minimum edge size $m$ if $2 <
q\leq (1-1/m)\Delta^{1/(m-1)}$.
\end{cor}
\begin{proof}
We cannot tell if there is \emph{any} colouring for $q$ in this range, so
there can be no \emph{fpras}.
\end{proof}
\begin{rem}
It is clearly a weakness that our lower bound for approximate counting is
based entirely on an NP-completeness result. However, we note that the same
situation pertains for  graph colouring, which has been the subject of more
intensive study.
\end{rem}
\section{Conclusions}
We have presented an approach to the analysis of path coupling with
stopping times which improves on the method of~\cite{HV04} in most
applications. Our method may itself permit further development.

We apply the method to independent sets and $q$-colourings in hypergraphs
with maximum degree $\Delta$ and minimum edge size $m$. In the case of
independent sets, there seems scope for improving the bound $m\geq 2\Delta
+1$, but anything better than $m\geq \Delta+o(\Delta)$ would seem to
require new methods. For colourings, there is probably little improvement
possible in our result $q >\Delta$ for $m\geq 4$, but many questions remain
for $m\leq \Delta$. For example, even the ergodicity of the Glauber (or any
other) dynamics is not clearly established. For the most interesting case,
$m=3$, the bound $q>1.65\Delta$ (for large $\Delta$) can almost certainly
be reduced, but substantial improvement may prove difficult.

Our \#P-completeness results seem best possible for both of the problems we
consider. On the other hand, our lower bounds for hardness of approximate
counting seem very weak in both cases, and are far from our upper bounds.
These lower bounds can probably be improved, but we have no plausible
conjecture as to what may be the truth.

\section*{Acknowledgments}

We are grateful to Tom Hayes for commenting on an earlier draft of this paper,
and to Mary Cryan for useful discussions at an early stage of this work.

\appendix

\section*{Appendices}

\section{Edge cover is \#P-complete}
\label{app1}
\begin{proof}
We prove this by reduction from counting independent sets, using methods
similar to Bubley and Dyer~\cite{BD97}, where this result was claimed
without proof. Let $\CG$ be a class of $3$-regular graphs for which
counting independent sets is \#P-complete~\cite[Theorem 3.1]{G00}. Let
$G=(V,E)$, with $I_j(G)$ independent sets of size $j$ ($j=0,1,\ldots,n$).
Form $G'$ by subdividing each edge $e\in E$ with a new vertex $u_e$. Let
$U=\set{u_e:e\in E}$. Let $N_i(G')$ be the number of edge sets in $G'$
which leave exactly $i$ vertices in $V$ uncovered, but no vertex in $U$. In
particular, $N_0(G')$ is the number of edge covers of $G'$, and we assume
an oracle computing this quantity. Observe that the uncovered vertices in
$G'$ must form an independent set in $G$. Then it follows, similarly
to~\cite{BD97}, that
\[ 2^{3(n-2j)}I_j(G)\ =\ \sum_{i=j}^n \binom{i}{j}N_i(G').\]
Thus, if we can determine the $N_i(G')$, we can determine the number of
independent sets of all sizes in $G$. Let $N_{ij}(G')$ be the number of
edge sets of $G'$ in which $i$ vertices in $V$ and $j$ in $U$ are
uncovered ($i=0,\ldots,n,\,j=0,\ldots,3n/2$). Then
$N_{i}(G')=N_{i0}(G')$. We attach a copy $K_m^v$ of $K_m$ to each vertex
$v\in V$ and a copy $K_k^u$ of $K_k$ to each vertex $u\in U$. Call the
resulting graph $G_{mk}$.  Let $M_m$ be the number of edge covers of $K_m$,
then $M_{m-1}$ is the number of edge sets in $K_m$ which leave a fixed
vertex uncovered. We can show by inclusion-exclusion that
\begin{equation*}
    M_m\ =\ \sum_{i=0}^{m}(-1)^i \binom{m}{i}
    \binpower{2}{m-i}.
\end{equation*}
(See~\cite{S04}.) It is easy to show that that $M_m/M_{m-1}$ is a rapidly
increasing sequence (in fact $M_m/M_{m-1}\approx 2^{m-1}$ for large $m$),
and hence has a different value for every value of $m$. We have
\begin{align}
N_{0}(G'_{mk})\ &=\ \sum_{i=0}^n \sum_{j=0}^{3n/2} M_m^{\,i}
(M_m+M_{m-1})^{n-i} M_k^{\,j} (M_k+M_{k-1})^{3n/2-j}
N_{ij}(G').\label{app1-eq0}\\
&=\ M_m^{\,n}M_k^{\,3n/2}\sum_{i=0}^n
\bigg(1+\frac{M_{m-1}}{M_m}\bigg)^{n-i}\ \sum_{j=0}^{3n/2}
\bigg(1+\frac{M_{k-1}}{M_k}\bigg)^{3n/2-j}N_{ij}(G').\notag
\end{align}
By choosing any $(n+1)$ values of $m$ and any $(3n/2 +1)$ values of $k$, we
can determine all the $N_{ij}(G')$ by interpolation, and hence all the
$N_i(G')$. From these, we can determine all the $I_j(G)$, and hence
$\sum_{j=1}^n I_j(G)$, the total number of independent sets in $G$.

Note that the minimum degree of $G'_{mk}$, $\min\{m,k\}-1$, can be made as
large as we wish.
\end{proof}
\section{Weak edge colouring is \#P-complete}
\label{app2}
\begin{proof}
We use the same notation and construction as in Appendix~\ref{app1}. Now
$\CG$ is a class of $3$-regular graphs for which vertex $q$-colouring is
\#P-complete~\cite[Theorem 2.2]{G00}. Let $N_i(G')$ denote the number of
edge colourings of $G'$ with $i$ monochromatic vertices, so $N_0(G')$ is
the number of weak edge colourings of $G'$, and we assume an oracle for
this. Let $N_{ij}(G')$ be the number edge colourings of $G'$ in which $i$
vertices in $V$ and $j$ in $U$ are monochromatic
($i=0,\ldots,n,\,j=0,\ldots,3n/2$). Now observe that $N_{n0}(G')$ is equal
to the number of proper \emph{vertex} $q$-colourings of $G$, $Q_q(G)$ say.
In every colouring counted in $N_{n0}(G')$, every vertex is monochromatic
and adjacent vertices receive different colours. Again we attach a copy
$K_m^v$ of $K_m$ to each vertex $v\in V$, and a copy $K_k^u$ of $K_k$ to
each vertex $u\in U$, to give $G_{mk}$. Let $M_m$ be the number of weak
colourings of $K_m$, and $M'_m$ the number of edge colourings of $K_m$ with
a given monochromatic vertex. Now we have
\begin{equation*}
    M_m\ =\ \binpower{q}{m}+q\sum_{i=1}^{m}(-1)^i \binom{m}{i}
    \binpower{q}{m-i},\qquad
    M'_m\ =\ q\sum_{i=0}^{m-1}(-1)^i \binom{m-1}{i}
    \binpower{q}{m-i-1}.
\end{equation*}
Again the sequence $M_m/M'_m$ increases rapidly ($M_m/M'_m\approx q^{m-2}$
for large $m$), and takes a different value for every $m$ when $q\geq 2$.
Now, as in (\ref{app1-eq0}),
\begin{equation*}
N_{0}(G'_{mk})\ =\ M_m^{\,n}M_k^{\,3n/2}\sum_{i=0}^n
\bigg(1+\frac{M'_m}{M_m}\bigg)^{n-i}\ \sum_{j=0}^{3n/2}
\bigg(1+\frac{M'_k}{M_k}\bigg)^{3n/2-j}N_{ij}(G').
\end{equation*}
Hence, choosing $(n+1)$ values of $m$ and $(3n/2 +1)$ values of $k$, we can
determine all the $N_{ij}(G')$ by interpolation. In particular, we can
determine $N_{n0}(G')=Q_q(G)$.

Again the minimum degree of $G'_{mk}$, $\min\{m,k\}-1$, can be made
arbitrarily large.
\end{proof}

\begin{thebibliography}{99}

\bibitem{BK03} P. Berman and M. Karpinski, Improved approximation lower bounds
on small occurrence optimization, \emph{Electronic Colloquium on
Computational Complexity} \textbf{10} (2003), Technical Report
TR03-008.

\bibitem{B01} R. Bubley,
\emph{Randomized algorithms: approximation, generation and counting},
Springer-Verlag, London, 2001.

\bibitem{BD97} R. Bubley and M. Dyer, Graph orientations with no sink
and an approximation for a hard case of \#SAT, in \emph{Proc.
8${}^\textrm{th}$ Annual ACM-SIAM Symposium on Discrete Algorithms
(SODA~1997)}, SIAM, 1997, pp.~248--257.

\bibitem{BDGJ99} R. Bubley, M. Dyer, C. Greenhill, and M. Jerrum, On approximately
counting colourings of small degree graphs, \emph{SIAM Journal on Computing}
\textbf{29} (1999), 387--400.

\bibitem{BG04} A. Bulatov and M. Grohe,
The complexity of partition functions, in \emph{Proc. 31st International Colloquium
on Automata, Languages and Programming (ICALP 2004)}, Springer, 2004,
pp.~294--306.

\bibitem{DGKR03}  I. Dinur, V. Guruswami, S. Khot and O. Regev,
A new multilayered PCP and the hardness of hypergraph vertex cover, in
\emph{Proc. 35${}^\textrm{th}$ ACM Symposium on Theory of Computing (STOC
2003)}, ACM, 2003, pp. 595--601.

\bibitem{DRS02}   I. Dinur, O. Regev and C. Smyth,
The hardness of {3-uniform} hypergraph coloring, in \emph{Proc.
43${}^\texttt{rd}$ Symposium on Foundations of Computer Science ({FOCS
2002})}, IEEE, 2002, pp. 33--42

\bibitem{DFHV04} M. Dyer, A. Frieze, T. Hayes and E. Vigoda, Randomly
coloring constant degree graphs, in \emph{Proc. 45${}^\textrm{th}$ Annual
IEEE Symposium on Foundations of Computer Science (FOCS 2004)}, IEEE, 2004,
pp.~582--589.

\bibitem{DFJ99} M. Dyer, A. Frieze and M. Jerrum, On counting independent
sets in sparse graphs, \emph{SIAM Journal on Computing} \textbf{31} (2002),
1527--1541.

\bibitem{DGGJM01} M. Dyer, L. Goldberg, C. Greenhill,
M. Jerrum and M. Mitzenmacher, An extension of path coupling and its
application to the Glauber dynamics for graph colorings,  \emph{SIAM
Journal on Computing} \textbf{30} (2001), 1962--1975.

\bibitem{DG98} M. Dyer and C. Greenhill, A more rapidly mixing Markov
chains for graph colouring, \emph{Random Structures and Algorithms}
\textbf{13} (1998), 210--217.

\bibitem{DG00a} M. Dyer and C. Greenhill, On Markov chains for independent
sets, \emph{Journal of Algorithms} \textbf{35} (2000), 17--49.

\bibitem{DG00b} M. Dyer and C. Greenhill, The complexity of counting graph
homomorphisms, \emph{Random Structures and Algorithms}, \textbf{17} (2000), 260--289.
See also Corrigendum, \emph{Random Structures and Algorithms} \textbf{25} (2004), 346--352.

\bibitem{EHK98} T. Emden-Weinert, S. Hougardy and B. Kreuter,
Uniquely colourable graphs and the hardness of colouring graphs of large girth,
\emph{Combinatorics, Probability and Computing} \textbf{7} (1998), 375--386.

\bibitem{GJ79} M. Garey and D. Johnson, \emph{Computer and intractability},
W. H. Freeman and Company, 1979.

\bibitem{G00} C. Greenhill, The complexity of counting colourings and
independent sets in sparse graphs and hypergraphs, \emph{Computational Complexity}
\textbf{9} (2000), 52--73.

\bibitem{HV04}
T. Hayes and E. Vigoda, Variable length path coupling, in \emph{Proc.
15${}^\textrm{th}$ Annual ACM-SIAM Symposium on Discrete Algorithms (SODA
2004)}, SIAM, 2004, pp.~103--110.

\bibitem{HL98} T. Hofmeister and H. Lefmann, Approximating maximum
independent sets in uniform hypergraphs, \emph{Proc. 23${}^\textrm{rd}$
International Symposium on Mathematical Foundations of Computer Science
(MFCS 1998)}, Lecture Notes in Computer Science \textbf{1450}, Springer,
1998, pp. 562--570.

\bibitem{JLR00} S. Janson, T. \L uczak and A. Ruci\'nski, \emph{Random graphs},
Wiley-Interscience, New York, 2000.

\bibitem{J95} M. Jerrum, A very simple algorithm for estimating the number of
$k$-colorings of a low-degree graph, \emph{Random Structure and Algorithms}
\textbf{7} (1995), 157--165.

\bibitem{J03} M. Jerrum, \emph{Counting, sampling and integrating: algorithms and
complexity}, ETH Z\"urich Lectures in Mathematics, Birkh\"auser, Basel,
2003.

\bibitem{KNS01} M. Krivelevich, R. Nathaniel and B. Sudakov,
Approximating coloring and maximum independent sets in 3-uniform
hypergraphs, in \emph{Proc. 12${}^\textrm{th}$ Annual ACM-SIAM Symposium on
Discrete Algorithms, (SODA~2001)}, SIAM, 2001, pp. 327--328.

\bibitem{LV99} M. Luby and E. Vigoda, Fast convergence of the Glauber dynamics
for sampling independent sets, \emph{Random Structures and
Algorithms} \textbf{15} (1999), 229--241.

\bibitem{M01} M. Molloy, Very rapidly mixing Markov chains for $2\Delta$-coloring
and for independent sets in a graph with maximum degree 4,
\emph{Random Structures and Algorithms} \textbf{18} (2001),
101--115.

\bibitem{MR01} M. Molloy and B. Reed, Colouring graphs when the
number of colours is nearly the maximum degree, in \emph{Proc.
33${}^\textrm{rd}$ Annual ACM Symposium on Theory of Computing (STOC
2001)}, ACM, 2001, pp. 462--470.

\bibitem{SS97} J. Salas and A. Sokal, Absence of phase transition
for anti-ferromagnetic Potts models via the Dobrushin uniqueness theorem,
\emph{Journal of Statistical Physics} \textbf{86} (1997), 551--579.

\bibitem{S04} N. Sloane, Sequence A006129,
\emph{The on-line encyclopedia of integer sequences},
 2004. Published at
\texttt{http://www.research.att.com/$\sim$njas/sequences/}.

\bibitem{V01} E. Vigoda, A note on the Glauber dynamics for sampling
independent sets, \emph{The Electronic Journal of Combinatorics}
\textbf{8}, R8(1), 2001.
\end{thebibliography}
\end{document}